\documentclass[11pt,draft]{elsart}
\usepackage{amssymb,amsmath}
\usepackage{mathrsfs}
\usepackage[english]{babel}
\usepackage{color}
\newtheorem{theo}{Theorem}[section]
\newtheorem{lema}{Lemma}[section]
\newtheorem{propt}{Proposition}[section]

\newcommand\XX[1]{\mathbb{#1}}

\newcommand\bd{\begin{definicion}{\bf }}
\newcommand\ed{\end{definicion}}
\newcommand\bl{\begin{lema}{\bf }}
\newcommand\el{\end{lema}}
\newcommand\bp{\begin{propt}{\bf }}
\newcommand\ep{\end{propt}}
\newcommand\bt{\begin{theo}{\bf }}
\newcommand\et{\end{theo}}
\newcommand\bdm{\begin{proof}}
\newcommand\edm{\end{proof}}
\newcommand\bn{\begin{nota}{\bf }}
\newcommand\en{\end{nota}}
\newcommand\bc{\begin{corolary}{\bf }}
\newcommand\ec{\end{corolary}}

\newtheorem{nota}{Remark}[section]

\newtheorem{definicion}{Definition}[section]
\newtheorem{corolary}{Corollary}[section]
\newenvironment{proof}{{\bf Proof:\ }}{\hfill
$\Box$}

\date{\today}

\journal{Journal of Approximation Theory}

\begin{document}

\begin{frontmatter}
\title{The semiclassical--Sobolev orthogonal polynomials:
a general approach}
\author{R.S. Costas--Santos}
\address{Department of Mathematics, University of California,
Santa Barbara, CA}
\author{J.J. Moreno--Balc\'{a}zar\corauthref{cor}}
\address{Departamento de Estad\'{i}stica y Matem\'{a}tica Aplicada,
\\ Universidad de Almer\'{i}a,  Spain \\
 Instituto Carlos I de F\'{\i}sica Te\'orica y Computacional}
\corauth[cor]{Corresponding author: e--mail: balcazar@ual.es}

\begin{abstract} We  say
that the polynomial sequence $(Q^{(\lambda)}_n)$ is a
se\-mi\-cla\-ssical Sobolev polynomial sequence when it is
orthogonal with respect to the inner product
$$
\left\langle p, r \right\rangle _S=\left\langle {{\bf u}}
,{p\, r}\right \rangle +\lambda \left \langle {{\bf u}},
{{\mathscr D}p \,{\mathscr D}r}\right\rangle,
$$
where ${\bf u}$ is a semiclassical linear functional, ${\mathscr
D}$ is the differential, the difference or the $q$--difference
operator, and $\lambda$ is a positive constant.

In this paper we get algebraic and differential/difference
properties for such polynomials as well as algebraic relations
between them and the polynomial sequence orthogonal with respect
to the semiclassical functional $\bf u$.

The main goal of this article is to give a general approach to the
study of the polynomials orthogonal with respect to the above
nonstandard inner product regardless of the type of operator
${\mathscr D}$ considered.  Finally, we illustrate our results
by applying them to some known families of Sobolev orthogonal
polynomials as well as to some new ones introduced  in this paper
for the first time.
\end{abstract}

\end{frontmatter}

\section{Introduction}
In the pioneering work \cite{lew1}, D.C. Lewis introduced the
inner products involving derivatives in order to obtain the least
squares approximation of a function and its derivatives. Since
that article, a large amount of papers have appeared on orthogonal
polynomials with respect to inner products involving derivatives
(and later differences or $q$--differences), the so--called
Sobolev (or $\Delta$--Sobolev or $q$--Sobolev) orthogonal
polynomials.

We can see the evolution of this theory through different surveys
from 1993 to 2006 (ordered by year: \cite{maalre}, \cite{mei1},
\cite{fink1}, \cite{fink2}, and \cite{mamo}). In the most recent
surveys some directions for the future are proposed.

The most studied cases correspond to the orthogonal polynomials
with respect to an inner product of the form
\begin{equation} \label{intr:ip}
\langle p, r \rangle_S= \langle
\mathbf{u}_0, p\, r \rangle
+\lambda \langle \mathbf{u}_1,
{\mathscr D}p \,{\mathscr D}r
\rangle, \qquad p, r \in \XX P,
\end{equation}
where $ \mathbf{u}_0$ and $\mathbf{u}_1$ are definite positive
linear functionals, and $\mathscr{D}$ is the differential
operator. Later, in \cite{argoma}  the discrete
case was considered  where ${\mathscr D}=\Delta$ was taken
as the forward finite difference operator, and in \cite{agm}
the $q$--difference operator ${\mathscr D}={\mathscr D}_q$
was considered.

As we can observe in the surveys mentioned above, and in the
references therein, all the results have been obtained by
considering a fixed operator, that is, the differential
operator in the continuous case, the forward operator in the
discrete case and the $q$--difference operator in the $q$--Hahn
case. Notice that a
great part of the results have been obtained from knowing an explicit
representation for the linear functionals $\mathbf{u}_0$ and
$\mathbf{u}_1$, for example, when $(\mathbf{u}_0, \mathbf{u}_1)$
is a coherent (or symmetrically coherent) pair of functionals. The
coherence is a case of special relevance from the numerical point
of view (see \cite{iskonosa}), although according to the
classification of the coherent pairs (see \cite{argoma},
\cite{agm}, and \cite{mei2} for each case), the coherent
functionals are close the classical ones; the complete
classification of the coherent pairs in the continuous, discrete
and $q$--Hahn cases can be found in \cite{mei2}, \cite{argoma} and
\cite{agm}, respectively. Furthermore, a lot of papers consider
the case $\mathbf{u}:=\mathbf{u}_0= \mathbf{u}_1$ when this
functional is explicitly known, for example, if $\bf u$ is a very
classical functional (Jacobi, Laguerre, Meixner, etc.) or if $\bf
u$ is semiclassical (Freud--type). For the corresponding Sobolev
orthogonal polynomials algebraic, differential/difference and
asymptotic properties are obtained as well as a relation with the
standard  polynomials orthogonal with respect to
 $\mathbf{u}$.

In this paper we do not work in the framework of coherence
(either $\Delta$--coherence or $q$--coherence) in the sense of
coherence worked with in \cite{mei2} (\cite{argoma} or \cite{agm}),
although some families of coherent pairs of this type are included
in our approach. In fact, we consider the case
$\mathbf{u}:=\mathbf{u}_0=\mathbf{u}_1$ in (\ref{intr:ip}) where
$\mathbf{u}$ is a semiclassical linear functional and, therefore,
our approach is not related to the coherence case.

We think that the theory of Sobolev orthogonal polynomials is a
big puzzle in which many pieces have been connected during the last
few decades but it lacks a  general approach. Taking this into
account, our main objective is to give a unified approach for the
theory of Sobolev orthogonal polynomials regardless of the
operator considered. More concretely, we consider a semiclassical
linear functional $\mathbf{u}$ and the nonstandard inner product
\begin{equation} \label{1:3}
\langle p, r \rangle_S= \langle
\mathbf{u}, p \, r \rangle+
\lambda \langle \mathbf{u},
{\mathscr D}p \, {\mathscr D}r
\rangle,\qquad p, r \in  \XX P,
\end{equation}
and we  obtain  differential/difference properties of the
orthogonal polynomials  with respect to (\ref{1:3}) in a unified
way independently of the type of operator ${\mathscr D}$, i.e.,
the results hold for the continuous, discrete, and $q$--Hahn
cases. We also find some relations between the semiclassical
Sobolev orthogonal polynomials with respect to the inner product
 (\ref{1:3}) and the semiclassical orthogonal
polynomials with respect to $\bf u$. To do this, we  use recent
results, obtained in \cite{coma1,khe}.

In the last section of the paper we illustrate our results,
applying them to some known families of Sobolev orthogonal
polynomials (Jacobi--Sobolev and $\Delta$--Meixner--Sobolev ones) as
well as to some  new ones introduced  in this paper for the first
time such as $q$--Fre\-ud--Sobolev ones and another family related to a
$1$--singular semiclassical functional.

In conclusion, we  show that it is possible to unify the
continuous, the discrete and the $q$--semiclassical Sobolev orthogonal
polynomials by using a suitable notation, which we believe
could be used in the future for further research.

\section{Basic definitions and notation. Semiclassical orthogonal
polynomials} As we said in the introduction, the Sobolev
polynomials, the $\Delta$--Sobolev polynomials and the
$q$--Sobolev polynomials have been considered in the literature
with different approaches. Our main idea  is to establish the main
algebraic properties of the polynomials which are orthogonal with
respect to the inner product $\left\langle \cdot, \cdot
\right\rangle _S$ defined in \eqref{1:3} in a general way where
the linear functional $\bf u$ satisfies a  distributional equation
with polynomial coefficients:
\begin{equation} \label{2:1} {\mathscr
D}(\phi{\bf u})=\psi {\bf u},
\end{equation}
where ${\mathscr D}$ is the differential or the difference or the
$q$--difference operator. Notice that a functional $\bf u$
satisfying (\ref{2:1}) with $\phi$ and $\psi$ polynomials, $\deg
\psi\ge 1$, is called a {\sf semiclassical} functional.

Since, depending on the case, the
polynomials are orthogonal with respect to different linear
functionals, we need to introduce a rigorous notation to embed all
of them.

\begin{enumerate}
\item {\it The Lattice.} For the general case we  use the variable
$z$. Notice that in the $q$--Hahn case we need to replace $z$ by
$x(s)= c_1(q)q^s+c_2(q)q^{-s}+c_3(q)$, where $s=0$, 1, 2, $\dots$,
with $q\in \XX C$, $|q|\ne 0, 1$, and $c_1$, $c_2$ and $c_3$ could
depend on $q$. However, in this paper we only consider  the
$q$--Hahn case, and so we work with the $q$--linear lattice
$$x(s):=c_1(q)q^s+c_3(q), $$ where $s=0$, 1, 2, $\dots,  $ and $q\in \XX C$, $|q|\ne 0, 1.$

\item {\it The operators.} For the general case we  use five basic
operators: ${\mathscr D}$, its dual ${\mathscr D}^*$, the identity
$\mathscr{I}$, the shift operator $\mathscr{E}^+$, and
$\mathscr{E}^-$. For a better reading of the paper we collect
these operators acting on a polynomial $r\in \mathbb{P}$ in the
following table:

\bigskip

\begin{center}
\begin{tabular}{|c|c|c|c|c|}
  \hline
   & $({\mathscr D}r)(x)$ & $({\mathscr D}^*r)(x)$ & $(\mathscr{E}^+r)(x)$ & $(\mathscr{E}^-r)(x)$ \\
  \hline
  Continuous case & $\frac{d r(x)}{dx}$ & $\frac{d r(x)}{dx}$ & $r(x)$ & $r(x)$ \\
  Discrete case  & $(\Delta r)(x)$ & $(\nabla r)(x)$ & $r(x+1)$ & $r(x-1)$ \\
  $q$--Hahn case & $({\mathscr
D}_{q} r)(x)$ & $({\mathscr D}_{1/q} r)(x)$& $r(qx)$ & $r(x/q)$ \\
  \hline
\end{tabular}
\end{center}

\bigskip

\noindent where the backward and forward difference operators are
defined as \begin{align*}
 (\Delta r)(x)\!:&=(\mathscr{E}^+-
\mathscr{I})(r)(x)=r(x+1)-r(x), \\ (\nabla r)(x)\!:&=(\mathscr{I}-
\mathscr{E}^-)(r)(x)=r(x)-r(x-1), \end{align*} and the
$q$--difference operator ${\mathscr D}_q$ is (see, for example,
\cite[Eq. (2.3)]{hah})
\begin{equation} \nonumber \label{2:2} ({\mathscr D}_q r)(x)=
\left\{ \begin{array}{ll} \displaystyle \frac{r(qx)-r(x)}
{(q-1)x}, & \quad x\ne 0,\\[5mm]
r'(0), & \quad x=0, \end{array}\right. \qquad r\in \XX P, \quad
q\in \XX C, \ |q|\ne 0,1.
\end{equation}

 \item {\it The
constants.} Throughout the paper we  use some constants: $q$, $[n]$,
and $[n]^*$. In the continuous and discrete cases we set  $q=1$,
and $[n]=[n]^*=n$, and in the $q$--Hahn case we  set $q\in \XX C$,
$|q|\ne 0, 1$, $[n]=(q^n-1)/(q-1)$, and
$[n]^*=(q^{-n}-1)/(q^{-1}-1)$.
\end{enumerate}

Let $\XX{P}$ be the linear space of polynomials and let $\XX{P}'$
be its algebraic dual space. We denote by $\left\langle
{{\bf{u}}}, p \right\rangle$ the duality bracket for ${\bf{u}} \in
\XX{P}'$ and $p \in \XX{P}$, and by $({\bf{u}})_n=\left\langle
{{\bf{u}}} ,{x^n}\right\rangle $, with $n \ge 0$, the canonical
moments of $\bf{u}$.

\bd For $\bf u \in \XX{P}'$,  $\pi \in \XX{P} $, and  $c \in
\mathbb{C}$, let $\pi\bf u$, $(x-c)^{-1}\bf u$, and ${\mathscr
D}^* \bf u$ be the linear functional defined by
$$
\left\langle {\pi\,{\bf{u}}},p \right\rangle:=\left\langle
{{\bf{u}}}, {\pi\,p}\right\rangle, \quad p \in \XX{P},
$$
\begin{equation*} \label{2:3}
\left\langle {(x-c)^{-1} {\bf{u}}},p \right\rangle:= \left\langle
{{\bf{u}}},{\frac{p(x)-p(c)} {x-c}}\right\rangle\, , \quad p \in
\XX{P},
\end{equation*}
\begin{equation} \label{2:4}
\left\langle {{\mathscr D}^*{\bf u}}, p\right\rangle:=
-q\left\langle {{\bf u}},{{\mathscr D} p}\right\rangle, \quad p\in
\XX P,
\end{equation}
and thus $\left\langle {{\mathscr D}{\bf u}}, p\right
\rangle=-q^{-1}\left\langle {{\bf u}} ,{{\mathscr D}^*
p}\right\rangle$.  \ed

Notice that for any ${\bf u}\in
\XX P'$,
\begin{equation*}\label{2:6}
(x-c)((x-c)^{-1} {\bf{u}})=\bf{u},
\end{equation*}
and
\begin{equation*}\label{2:5}
(x-c)^{-1}((x-c){\bf{u}})={\bf{u}}
-({\bf{u}})_0\delta_c,
\end{equation*}
where $\delta_c$  is the Dirac delta functional defined by
$\langle \delta_c,p\rangle:=p(c)$,  $p\in \mathbb P.$

Taking these definitions into account  we get for any ${\bf u}\in
\XX P'$ and any polynomial $\pi$
\begin{align}
{\mathscr D}(\pi {\bf u})=& (\mathscr{E}^+\pi){\mathscr D}{\bf
u}+({\mathscr D}\pi)
{\bf u}, \label{2:7}\\
{\mathscr D}^*(\pi{\bf u})=& (\mathscr{E}^-\pi){\mathscr D}^*{\bf
u}+({\mathscr D}^*\pi) {\bf u}. \label{2:8}
\end{align}
Furthermore, for any two polynomials, $\pi$ and $\xi$, the
following relations are fulfilled:
\begin{align}
\label{2:9} \displaystyle {\mathscr D}(\pi
\xi)&=
(\mathscr{E}^+ \pi){\mathscr D}
\xi+({\mathscr D}\pi) \xi,\\
\label{2:10} \displaystyle {\mathscr D}^* (\pi \xi)&=
(\mathscr{E}^-\pi){\mathscr D}^*\xi+({\mathscr D}^* \pi) \xi.
\end{align}

\bd Given an integer $\delta$, a linear functional $\bf u$, and a
family of polynomials $(p_n)$, we say that the polynomial sequence
$(p_n)$ is {\sf quasi--orthogonal of order $\delta$ with respect
to} $\bf u$ when the following properties are fulfilled:
\begin{itemize} \item If $|n-m|\ge \delta+1$, then $\left\langle
{{\bf u}}, {p_np_m} \right\rangle=0$. \item If $|n-m|=\delta$,
then $\left\langle {{\bf u}}, {p_np_m} \right\rangle\ne 0$.
\end{itemize} \ed

Notice that if $\bf u$ is a semiclassical functional and $(p_n)$
is orthogonal with respect to $\bf u$, i.e. $\delta=0$, then
$(p_n)$ is called {\sf semiclassical orthogonal polynomial
sequence}.

Now we  introduce some results regarding semiclassical orthogonal
polynomials which allow us to obtain several results concerning
 Sobolev orthogonal polynomials in Sections 3 and 4. There is a
lot of literature on semiclassical polynomials and an
authoritative paper on this topic is \cite{maroni-1999}.

First of all,  we define the concept of an admissible pair which
appears frequently linked with the concept of the semiclassical linear
functional.

\bd \label{def3.2} We say that the pair of polynomials $(\phi,
\psi)$ is an {\sf admissible pair} if one of the following
conditions is satisfied: \begin{itemize} \item  $\deg \psi\ne \deg
\phi-1$. \item  $\deg \psi=\deg \phi-1$, with $$
a_p+q^{-1}[n]^*b_t\ne 0,\qquad  n\ge 0, $$ where $a_p$ and $b_t$
are the leading coefficients of $\psi$ and $\phi$, respectively.
\end{itemize}
We  denote by $\sigma:= \max\{\deg\phi-2,\deg \psi -1\}$ the {\sf
order of the linear functional} $\bf u$ with respect to the
admissible pair $(\phi,\psi)$. The {\sf class} of such functional
is the minimum of the order from among all the admissible pairs
(see \cite{khe,mar7}). \ed \bn \label{rem3.1} Taking into account
the two possible situations for the admissibility we use the
following notation to unify them. We  write
$$
\phi(z)=b_{\sigma+2}z^{\sigma+2}+
\cdots, \qquad \psi(z)=a_{\sigma+
1}z^{\sigma+1}+\cdots,
$$
where:
\begin{itemize}
\item If $t=\sigma+ 2$ and $p<\sigma+1$, then we set
$b_{\sigma+2}=b_t$ and $a_{\sigma+1}=a_{\sigma}=\cdots=a_{p+1}=0$.
\item If $p=\sigma+ 1$ and $t<\sigma+2$, then we  set
$a_{\sigma+1}=a_p$ and
$b_{\sigma+2}=b_{\sigma+1}=\cdots=b_{t+1}=0$. \end{itemize} \en
\bd \label{def3.3} We say that the distributional equation
\eqref{2:1} has an $n_0$--{\sf singularity} if there exists a
non--negative integer $n_0$ such that $a_{\sigma+1}+ q^{-1}[n_0]^*
b_{\sigma+2}=0$. Otherwise we say that the distributional equation
is {\sf regular}. \ed The following result will be very useful for
our purposes. In the next result  and in Lemma \ref{lem3.1.1} we
give some details in their proofs for each case, namely the continuous,
discrete and $q$--Hahn cases, with the objective of making them more
readable. It is also important to point out that all the results
presented throughout the paper are valid for all the cases with no
additional restrictions.

\bt \label{the3.0} Let $\bf u$ be a linear functional
satisfying the distributional equation
$$
{\mathscr D}(\phi{\bf u})=
\psi {\bf u},
$$
where $\phi$ and $\psi$ are polynomials; then $\bf u$ also
fulfills the distributional equation
\begin{equation}  \nonumber\label{3:5.1}
{\mathscr D}^*(\widetilde {\phi}{\bf u})= \widetilde {\psi}{\bf
u}, \end{equation} where $\widetilde \phi$ and $\widetilde \psi$
are polynomials, with $\deg \widetilde {\phi}\le \sigma+2$, and
$\deg \widetilde {\psi}\le \sigma+1$, with $\sigma=\max\{\deg
\phi-2,\deg\psi-1\}$. \et \bdm The  continuous case is trivial
since ${\mathscr D}={\mathscr D}^*=d/dx.$

In the discrete case ${\mathscr D}=\Delta$ and ${\mathscr D}^*=
\nabla$, so
$$
{\mathscr D}(\phi{\bf u})=\psi {\bf u} \iff
\Delta(\phi {\bf u})={\mathscr E}^+(\phi {\bf u})-\phi
{\bf u}=\psi {\bf u} \iff
{\mathscr E}^+(\phi {\bf u})=(\psi+\phi){\bf u}.
$$
Observe that, due to the property of the shift operator ${\mathscr
E}^+$, the latter expression is also equivalent to $\phi {\bf
u}={\mathscr E}^-((\psi+\phi){\bf u}). $ Thus, setting $\widetilde
\phi:=\phi+\psi$ and $\widetilde \psi=\psi$, we get
$$
{\mathscr D}^*(\widetilde \phi {\bf u})=\nabla(
\widetilde \phi {\bf u})=(\phi+\psi){\bf u}-
{\mathscr E}^-((\phi+\psi){\bf u})=\psi {\bf u}
=\widetilde \psi {\bf u}.
$$
The $q$--Hahn case is a little more complicated, due to the
notation, but the proof is almost identical.
\\
In this case ${\mathscr D}=\Delta/\Delta x(s)$, and
${\mathscr D}^*=\nabla/\nabla x(s)$;
thus the  distributional equation for ${\mathscr D}$
can be written as
$$
\frac{{\mathscr E}^+(\phi {\bf u})}{\Delta x}- \frac{\phi {\bf
u}}{\Delta x}= \frac{(\Delta x) \psi {\bf u}}{\Delta x} \quad \iff
\quad \frac{{\mathscr E}^+(\phi {\bf u})}{\Delta x}=
\frac{(\phi+(\Delta x)\psi)}{\Delta x}{\bf u}.
$$
Then,
$$\begin{array}{rl}
\dfrac{({\mathscr E}^-(\phi+(\Delta x) \psi))}
{\nabla x}\nabla {\bf u} =& \dfrac{({\mathscr E}^-(
\phi+(\Delta x) \psi))}{\nabla x}{\bf u}-\dfrac{(
{\mathscr E}^-(\phi+(\Delta x) \psi))}
{\nabla x}{\mathscr E}^-{\bf u}\\[3mm]
& = \displaystyle \dfrac{({\mathscr E}^-(\phi+(\Delta x)
\psi))} {\nabla x}{\bf u}-\frac{\phi {\bf u}}{\nabla x} =
\left({\mathscr E}^-\psi-\frac{\nabla \phi} {\nabla x}\right)
{\bf u}. \end{array}$$
Thus, using the above expression, \eqref{2:8}, and setting
$\widetilde \phi:=\phi+(\Delta x)\psi$, we get
$$
{\mathscr D}^*(\widetilde \phi {\bf u})=\dfrac{\nabla\Big(
(\phi+(\Delta x)\psi){\bf u}\Big)} {\nabla x} =\left({\mathscr
E}^-\psi-\frac{\nabla \phi} {\nabla x}\right){\bf u}+\dfrac{\nabla
(\phi+(\Delta x)\psi)}{\nabla x}{\bf u} =q\psi {\bf u}=\widetilde
{\psi}{\bf u}.
$$
Therefore, it is a straightforward computation to check that both
$\widetilde \phi$ and $\widetilde \psi$ are polynomials, with
$\deg \widetilde {\phi}\le \sigma+2$, and $\deg \widetilde
{\psi}\le \sigma+1$ where $\sigma=\max\{\deg \phi-2,\deg\psi-1\}$.
\edm

\bn \begin{itemize} \item Although the reciprocal of this result
is also true, and it can be proved in a straightforward way, we have
decided not to include it since it is not used in this paper.

\item In the sequel, given a pair $(\phi,\psi)$ we use the notation
$\widetilde \phi:=\phi+\psi$ and  $\widetilde \psi:=\psi$ for
the discrete case, and  $\widetilde \phi:=\phi+(\Delta x)\psi$
and $\widetilde \psi:=q\psi$ for the $q$--Hahn case.
\end{itemize} \en

The inner product associated
with the semiclassical
functional $\bf u$,
$$
\left\langle p, r \right\rangle:=\left\langle {{\bf u}},
{p \, r}\right\rangle,
$$
is Hankel, i.e.,
\begin{equation} \label{omi} \langle {zp},r
\rangle= \langle p, {zr}\rangle,
\end{equation}
and then the corresponding monic semiclassical orthogonal polynomial
sequence $(p_n)$ satisfies the recurrence relation
\begin{equation}  \label{3:5} zp_n(z)=p_{n+1}(z)+B_n
p_n(z)+C_np_{n-1}(z),\quad n\ge 0, \end{equation} with initial
conditions $p_{0}=1$, $p_{-1}=0$. Moreover the square of the norm
of $p_n$, $d^2_n=\left\langle {{\bf u}}, {p_n^2}\right\rangle$,
satisfies the relation
\begin{equation} \label{3:6}
d^2_n=C_n d^2_{n-1},\qquad n\ge 1.
\end{equation}

It is well known that there are several ways to characterize
semiclassical orthogonal polynomials.  In our framework, we  use a
characterization of these polynomials that  is usually called the
first structure relation. \bt \label{lem3.1} \cite{mar} \&
\cite[Prop. 3.2]{khe}  Let $(p_n)$ be a sequence of monic
orthogonal polynomials with respect to $\bf u$, and $\phi$ a
polynomial of degree $t$. Then,  the following statements are
equivalent:
\begin{enumerate} \item[\bf (i)] There exists a non-negative
integer $\sigma$ such that
\begin{align*}
\phi(z){\mathscr D}p_{n+1}(z)&=\displaystyle \sum_{\nu=
n-\sigma}^{n+t}\lambda_{n,\nu} p_\nu(z),\quad n\ge \sigma,\\
\lambda_{n,n-\sigma}&\ne 0, \quad n\ge\sigma+1.
\end{align*}
\item[\bf (ii)] There exists a polynomial $\psi$, with $\deg
\psi=p\ge 1$, such that
\begin{equation} \nonumber \label{3:9} {\mathscr
D}(\phi{\bf u})=\psi {\bf u},
\end{equation}
where the pair $(\phi,\psi)$ is an admissible pair.
\end{enumerate} \et

We also use  a recent result that  establishes a nice relation
between the semiclassical orthogonal polynomials $(p_n)$ and the
polynomials $({\mathscr D}p_{n+1}).$

 \bt \label{the3.1} \cite[Th. 3.2]{coma1}
\& \cite[Th. 3.5]{sfma} Let $(p_n)$ be a sequence of monic
orthogonal polynomials with respect to $\bf u$, and $\phi$ a
polynomial of degree $t$. Then,  the following statements are
equivalent:
\begin{enumerate}
\item[\bf (i)] There exist three non--negative integers, $\sigma$,
$p$, and $r$, with $p\ge 1$, $r\ge \sigma+t+1$, and
$\sigma=\max\{t-2,p-1\}$, such that \begin{equation} \nonumber
\label{3:10} \sum_{\nu=n- \sigma}^{n+\sigma} \xi_{n,\nu} p_\nu(z)=
\sum_{\nu=n-t}^{n+\sigma} \varsigma_{n,\nu}p_\nu^{[1]}(z),
\end{equation} where $p_n^{[1]}(z):=[n+1]^{-1}({\mathscr D}
p_{n+1}) (z)$, $\xi_{n,n+\sigma}=\varsigma_{n,n+\sigma}=1$, $n\ge
\max\{\sigma,t+1\}$, $\xi_{r,r-\sigma} \varsigma_{r,r-t}\ne 0$, $$
\langle {\mathscr D}(\phi {\bf u}), p_m\rangle=0, \quad p+1 \le
m\le 2 \sigma+t+1, \quad \langle {\mathscr D}(\phi {\bf u}),
p_p\rangle\ne 0,
$$
and if $p=t-1$, then $ \lim_{q\uparrow 1} \langle {\bf
u},p_p^2\rangle^{-1} \langle {\bf u},\phi {\mathscr D} p_p\rangle
\ne m b_t$, $m\in \XX N^*,$ and where $b_t$ is the leading
coefficient of $\phi$ (admissibility condition). \item[\bf (ii)]
There exists a polynomial $\psi$, with $\deg \psi=p\ge 1$, such
that
\begin{equation} \nonumber \label{3:11}
{\mathscr D}(\phi{\bf u})=\psi {\bf u},
\end{equation}
where the pair $(\phi,\psi)$ is
an admissible pair. \end{enumerate} \et

In this paper we use some linear functionals related to the
semiclassical functional $\bf u$. Now, we introduce them and give
some properties that they satisfy.

\bl \label{lem3.1.1} Let $\bf u$ be a linear functional satisfying
\eqref{2:1}; then the linear functional ${\mathscr E}^+ (\phi {\bf
u})$ satisfies the following distributional equation:
\begin{equation*} \label{eq:dis11}
{\mathscr D}(\phi\, {\mathscr E}^+ (\phi {\bf u}))=\big({\mathscr
D}\phi+q{\mathscr E}^+\psi\big) {\mathscr E}^+ (\phi {\bf u}),
\end{equation*}
setting $q=1$ in both cases, continuous and discrete. \el

\bdm The result follows by applying \eqref{2:7}.
In fact, in the continuous case, since ${\mathscr E}^+$ is the
identity, we get
$${\mathscr D}(\phi\, {\mathscr E}^+ (\phi {\bf u}))
=\frac{d}{dx}(\phi \phi {\bf u})
=\frac{d}{dx}(\phi) \phi {\bf u}+
\phi \frac{d}{dx}(\phi {\bf u})=
\big({\mathscr D}\phi+{\mathscr E}^+\psi\big) {\mathscr E}^+
(\phi {\bf u}).
$$
In the discrete case, ${\mathscr E}^+$ is the shift operator, so
applying \eqref{2:7} we get
$$\begin{array}{rl}
{\mathscr D}(\phi\, {\mathscr E}^+ (\phi {\bf u}))=&
\Delta(\phi\, {\mathscr E}^+ (\phi {\bf u})) =
{\mathscr E}^+ \phi\, \Delta({\mathscr E}^+ (\phi
{\bf u}))+\Delta \phi\, {\mathscr E}^+ (\phi {\bf u})
\\ =& {\mathscr E}^+\phi\, {\mathscr E}^+ (\psi {\bf
u})+\Delta \phi\, {\mathscr E}^+ (\phi {\bf u})\\
=& \left({\mathscr E}^+\psi+\Delta \phi \right)
{\mathscr E}^+ (\phi{\bf u}).
\end{array}$$
In the $q$--Hahn case, ${\mathscr E}^+$ is also  the shift
operator. Thus applying \eqref{2:7} we get
$${\mathscr D}(\phi\, {\mathscr E}^+ (\phi {\bf u}))
=\frac{\Delta(\phi\, {\mathscr E}^+ (\phi {\bf u}))}{\Delta x(s)}
=\frac 1{\Delta x(s)}\big(\Delta \phi\, {\mathscr E}^+ (\phi {\bf
u})+ {\mathscr E}^+\phi\, \Delta({\mathscr E}^+ (\phi {\bf
u}))\big).
$$
Taking into account that in this case ${\mathscr D}{\mathscr
E}^+=q{\mathscr E}^+{\mathscr D}$, since the lattice $x(s)$ is
$q$--linear, we get
$$
{\mathscr D}(\phi\, {\mathscr E}^+ (\phi {\bf u}))= \big({\mathscr
D}\phi\, {\mathscr E}^+(\phi{\bf u})+ q{\mathscr E}^+\phi
\,{\mathscr E}^+{\mathscr D}(\phi{\bf u})\big) =\big({\mathscr
D}\phi+q{\mathscr E}^+\psi \big){\mathscr E}^+ (\phi{\bf u}).
$$
\edm

\bl \label{lem3.1.2} Let $\bf u$ be a linear functional satisfying
\eqref{2:1}, and let  $\bf v$ be a linear functional such that
${\bf u}= {\mathscr E}^+(\phi {\bf v})$; then $\bf v$ satisfies an
analogous distributional equation. \el

\bdm It follows from making computations in  \eqref{2:1} where it
is necessary to take into account how the operators ${\mathscr D}$
and ${\mathscr E}^+$ act, that is, ${\mathscr D}{\mathscr
E}^+=q{\mathscr E}^+{\mathscr D}.$ In fact, we obtain
$$
{\mathscr D}(q({\mathscr E}^-\phi) \phi \, {\bf v}) =({\mathscr
E}^-\psi)\phi \, {\bf v}.
$$
\edm

So, the previous  lemmas provide motivation for defining a family of
linear functionals which will be used throughout the paper. \bd
\label{def3.4}  Given a semiclassical functional $\bf u$
satisfying \eqref{2:1}, for  $k\in \mathbb{Z}$ we define the
linear functional ${\bf u}_{k}$ in a constructive way as
\begin{equation*} \label{3:13}
\begin{array}{rl} {\bf u}_{k}=& \mathscr{E}^+(\phi_{k-1}
{\bf u}_{k-1}), \\
{\bf u}_{0}=& {\bf u}, \quad \phi_0=\phi,
\end{array}\end{equation*}
where $\phi_k$ is a multiple of $\phi_{k-1}.$
\ed  \bn \label{rem3.3}
 Notice that, by Lemma \ref{lem3.1.1} and Lemma \ref{lem3.1.2}, when
$\bf u$ is a semiclassical functional then ${\bf u}_k$ is a
semiclassical functional for every integer $k$.  On the other hand,
using Lemma 3.3 in \cite{mar7} for the continuous case and Lemma 2.2
in \cite{khe} for  discrete and $q$--Hahn cases,  the statement of
$\mathbf{u}_k$ in terms of $\mathbf{u}_{k-1}$ is well defined since
if $\mathbf{u}$ is any semiclassical functional satisfying the
distributional equations
$$
{\mathscr D}(\hat{\phi}_1 \mathbf{u})=\hat{\psi}_1 \mathbf{u}
\quad \mathrm{and} \quad {\mathscr D}(\hat{\phi}_2
\mathbf{u})=\hat{\psi}_2 \mathbf{u},$$
then for $\hat{\phi}:=\gcd(\hat{\phi}_1,\hat{\phi}_2)$ there exists a
polynomial $\hat{\psi}$ with $\deg \hat \psi\ge 1$ such that
${\mathscr D}(\hat{\phi} \mathbf{u})=\hat{\psi} \mathbf{u}.$ Thus,
observe that $\phi_k$ is not determined in a unique way. \en

As a consequence of Lemmas \ref{lem3.1.1} and \ref{lem3.1.2} and
some straightforward computations we can deduce:

\bt \label{lem3.2} Say we are given a linear functional $\bf u$ satisfying
\eqref{2:1}, that is,
$${\mathscr D}(\phi{\bf u})=\psi {\bf
u}, $$ where $\psi$ and $\phi$ are polynomials in the conditions
of Definition \ref{def3.2}. Then, for any integer $k$ and any
polynomial $\pi$ we get
\begin{equation} \nonumber \label{3:14}
{\mathscr D}(\pi \phi_k {\bf u}_k)=\widetilde \pi {\bf u}_k,
\end{equation}
where $\widetilde \pi$ is a polynomial of degree, at most, $\deg
\pi +\sigma_k+1$, $\sigma_k$ being  the order of the linear
functional ${\bf u}_k$ with respect to the pair $(\phi_k,\psi_k)$.
Furthermore, if the pair $(\phi_k,\psi_k)$ is admissible, then
$\deg \widetilde \pi=\deg \pi+ \sigma_k+1$. \et

Finally, we obtain a useful result for our purposes.

\bp \label{cor 3:2} If $(p_n)$ is a sequence of semiclassical
polynomials orthogonal with respect to the linear functional $\bf
u$ of order $\sigma$, then the sequence of polynomials $({\mathscr
D}p_{n+1})$ is quasi--orthogonal of order $\sigma$ with respect to
the semiclassical functional ${\bf u}_1$. \ep

\bdm Using (\ref{2:1}), (\ref{2:4}) and (\ref{2:10}),  we have
\begin{align*}
\langle {\bf u}_1, {\mathscr D}p_{n+1}x^m \rangle & = \langle
{\mathscr E}^+(\phi{\bf u}), {\mathscr D}p_{n+1}x^m \rangle=
\langle \phi {\bf u}, {\mathscr E}^-({\mathscr D}p_{n+1}x^m
)\rangle
\\ & =\langle \phi {\bf u}, {\mathscr D}^*(p_{n+1}x^m)\rangle-
\langle \phi {\bf u}, p_{n+1}{\mathscr D}^*x^m\rangle \\
&=-\langle {\bf u},p_{n+1}(q\psi x^m+\phi{\mathscr D}^*x^m)
\rangle=0,
\end{align*}
when $m+\sigma<n.$ \edm

\bn As far as we know, the reciprocal of Proposition 2.1 is no
longer true without the assumption that $\bf v$ is quasi-definite
where $\bf v$ is the operator with respect to which $({\mathscr
D}p_{n+1})$ is quasi--orthogonal. \en

For every integer $k$ we denote by $(p_n^{\{k\}})$  the sequence
of monic orthogonal polynomials  with respect to the linear
functional ${\bf u}_k$. Thus, using Theorem \ref{lem3.2} and
Proposition \ref{cor 3:2}, we have that the ${\mathscr D}
p^{\{k\}}_{n}$ are quasi--orthogonal of order $\sigma_k$ with
respect to ${\bf u}_{k+1}$. Therefore,  we get  for every integer
$k$,
\begin{equation} \label{3:17}
{\mathscr D} p_{n+1}^{\{k\}}
=\sum_{\nu=n-\sigma_{k}}^n \alpha_{n,k,\nu} p_\nu^{\{k+1\}},
\qquad k\in \XX Z.
\end{equation}

Then, using \eqref{3:17} and Theorem \ref{the3.1}, we deduce for
$k=-1$ the following  relation:
\begin{equation} \label{3:18}
\sum_{\nu=n-\sigma_{-1}}^{n+\sigma_{-1}} \xi^*_{n,\nu}
p^{\{-1\}}_\nu(z)=\sum_{\nu=n-t_{-1}-\sigma_{-1}}^{n+\sigma_{-1}}
\varsigma^*_{n,\nu}p_\nu(z), \end{equation} where
$\xi^*_{n,n+\sigma_{-1}-1}=\varsigma^*_{n, n+\sigma_{-1}}=1$,
$t_{-1}:=\deg \phi_{-1}$, and $\sigma_{-1}$ is the order of the
functional ${\bf u}_{-1}$ associated with the admissible pair
$(\phi_{-1},\psi_{-1})$.

Observe that in the classical case, that is, $\sigma=0$, using
Theorem \ref{lem3.2} we have  $\sigma_{-1}=0$.


\section{The semiclassical Sobolev orthogonal polynomials}

Let us consider the ${\mathscr D}$--Sobolev inner product defined
on $\XX P\times \XX P$ by
\begin{equation} \label{4:1}
\langle p, r \rangle_S=\langle {{\bf u}},
{p\, r}\rangle+\lambda \langle {{\bf u}}, {{\mathscr D} p\,
{\mathscr D} r}\rangle,
\end{equation}
where ${\bf u}$ is a semiclassical functional of order $\sigma$
and $\lambda \ge 0$. We denote by  $(p_n)$ the monic orthogonal
polynomial sequence with respect to ${\bf u}$ and by
$(Q_n^{(\lambda)})$ the OPS associated with the ($\mathscr
D$--Sobolev) inner product $\langle \cdot, \cdot \rangle_S$ which
we call the semiclassical Sobolev orthogonal polynomial sequence.

\bp \label{prop4.1} Let $(p_n)$ be a monic semiclassical
orthogonal polynomial sequence and let $(Q_n^{(\lambda)})$ be the
semiclassical Sobolev orthogonal polynomial sequence. The
following relation holds: \begin{equation} \label{4:2}
\sum_{\nu=n-\sigma_{-1}}^{n+\sigma_{-1}} \xi^*_{n,\nu}
p^{\{-1\}}_\nu(z)= Q^{(\lambda)}_{n+\sigma_{-1}}(z)+
\sum_{\nu=n-\sigma_{-1}-H^*}^{n-1+\sigma_{-1}} \theta_{n,\nu}\,
Q^{(\lambda)}_\nu(z), \quad n \ge \sigma_{-1}+H^*,
\end{equation}
where
$H^*:=\max\{t_{-1},\sigma_{-1}\},$
\begin{equation} \label{4:3}
\begin{split}{\bf d}^2_{n-\sigma_{-1}-H^*}
\theta_{n,n-\sigma_{-1}-H^*}
(\lambda)=& \varsigma^*_{n,n-t_{-1}-\sigma_{-1}}
d^2_{n-t_{-1}-\sigma_{-1}}
\delta_{H^*,t_{-1}}\\  +&\lambda [n-2\sigma_{-1}]
\hat\xi^*_{n,n-2\sigma_{-1}}d^2_{n-2\sigma_{-1}-1}
\delta_{H^*,\sigma_{-1}},
\end{split}\end{equation}
with ${\bf d}^2_n:=\left\langle {Q^{(\lambda)}_n(z)}
{Q^{(\lambda)}_n(z)} \right\rangle _S$, and the other coefficients
$\theta_{n,\nu}$ in (\ref{4:2}) can be computed recursively. \ep
\bdm If we apply \eqref{3:18} and expand it we get
$$
\sum_{\nu=n-\sigma_{-1}}^{n+\sigma_{-1}} \xi^*_{n,\nu}
p^{\{-1\}}_\nu(z)= \sum_{\nu=n-t_{-1}-\sigma_{-1}}^{n+\sigma_{-1}}
\varsigma^*_{n,\nu}p_\nu(z)
=Q^{(\lambda)}_{n+\sigma_{-1}}(z)+\sum_{i=0}^{n+\sigma_{-1}-1}
\theta_{n,i}Q^{(\lambda)}_i(z),
$$
where
$$
\theta_{n,i}=\frac{\left\langle {\displaystyle \sum_{\nu=n-
\sigma_{-1}}^{n+\sigma_{-1}}\xi^*_{n,\nu}
p^{\{-1\}}_\nu(z)},{Q^{(\lambda)}_i(z)}\right\rangle_S}
{\left\langle {Q^{(\lambda)}_i(z)},{Q^{(\lambda)}_i(z)}\right
\rangle_S}=\frac{\displaystyle \left\langle {\sum_{\nu=n-t_{-1}-
\sigma_{-1}}^{n+\sigma_{-1}} \varsigma^*_{n,\nu} p_\nu(z)},
{Q^{(\lambda)}_i(z)}\right\rangle_S}{{\bf d}^2_i}.
$$
Thus, using \eqref{3:17} and the orthogonality  property we get
for $0\le i<n-\sigma_{-1}-H^*$
$$\begin{array}{rl}
{\bf d}^2_i \,\theta_{n,i}= & \displaystyle \left\langle
{{\bf u}}, {\left(\sum_{\nu=n-t_{-1}-\sigma_{-1}}^{n+
\sigma_{-1}} \varsigma^*_{n,\nu}p_\nu(z)\right)
Q^{(\lambda)}_i(z)}\right\rangle\\[5mm] & +\displaystyle
\lambda \left\langle {{\bf u}}, {{\mathscr D} \left(
\sum_{\nu=n-\sigma_{-1}}^{n+\sigma_{-1}} \xi^*_{n,\nu}
p^{\{-1\}}_\nu(z) \right)\, ({\mathscr D} Q^{(\lambda)}_i)
(z)}\right\rangle \\[4mm]= & \displaystyle \lambda \left
\langle {{\bf u}} ,{\left(\sum_{\nu=n- 2\sigma_{-1}}^{n+
\sigma_{-1}}\hat \xi^*_{n,\nu} p_{\nu-1}(z)\right)\,
({\mathscr D} Q^{(\lambda)}_i)(z)}\right\rangle = 0,
\end{array}
$$
and therefore (\ref{4:2}) follows. To obtain (\ref{4:3}) we
only need to take $i=n-\sigma_{-1}-H^*$ in the above expression.

Finally,  if we apply this process for $n-\sigma_{-1}-H^*\le j\le
n+\sigma_{-1}-1$ we get an upper triangular linear system such
that we can compute the coefficients $\theta_{n,j}$ in a recursive
way and hence the result holds. \edm

Next, we particularize the above result for the classical case,
that is, we name as the classical Sobolev orthogonal polynomial
sequence the semiclassical Sobolev orthogonal polynomial
sequence where the linear functional involved has order
$\sigma=0,$ i.e., it is a classical functional. \bc \label{cor4.1}
Let $(p_n)$ be a classical orthogonal polynomial sequence and let
$(Q_n^{(\lambda)})$ be the classical Sobolev orthogonal
polynomial sequence. The following relation holds:
\begin{equation} \label{4:4} p_{n}^{\{-1\}}(z)=Q^{(\lambda)}_n(z)+
f_{n-1}(\lambda)Q^{(\lambda)}_{n-1}
(z)+e_{n-2}(\lambda)Q^{(\lambda)}_{n-2}(z), \quad n\ge 2,
\end{equation}
with $Q^{(\lambda)}_i(z)=p_i(z)$, $i=0,1, $ and
\begin{eqnarray} \label{4:5}
 e_{n-2}(\lambda)&=& \tilde
\epsilon_n \frac{d^2_{n-2}}{{\bf d}^2_{n-2}}, \quad
e_0(\lambda)=\tilde
\epsilon_2, \\[2mm] \nonumber \label{4:6}
f_{n-1}(\lambda)&=&\frac 1 {{\bf d}^2_{n-1}} \left\{\tilde
\delta_n d^2_{n-1}+\tilde \epsilon_n (\tilde \delta_{n-1}
-f_{n-2}(\lambda))d^2_{n-2}\right\}, \quad f_0(\lambda)=\tilde
\delta_1,
\end{eqnarray}
where
\begin{equation} \label{4:7}
p_{n}^{\{-1\}}(z)=p_n(z)+\tilde \delta_n
p_{n-1}(z)+\tilde \epsilon_n p_{n-2}(z).
\end{equation}
\ec

Furthermore, in the classical case we can give nonlinear
recurrence relations for the square of the norms of the Sobolev
polynomials and for the coefficients $e_n(\lambda)$ appearing in
\eqref{4:4}.

\bc \label{prop4.2} For $n\ge 2$, we have
\begin{equation} \label{4:8}
{\bf d}^2_n=d^2_n\!+\big(\!\lambda
[n]^2+\tilde \delta_n(\tilde \delta_n
-f_{n-1}(\lambda))\big)d^2_{n-1}
+\tilde \epsilon_n (\tilde \epsilon_n
- e_{n-2}(\lambda)-f_{n-1}(\lambda)
(\tilde \delta_{n-1}-f_{n-2} (\lambda)))
d^2_{n-2},
\end{equation}
with ${\bf d}^2_0=d^2_0$, and ${\bf d}^2_1=d^2_1+\lambda d^2_0$,
and
\begin{equation} \label{4:9}
e_n(\lambda)=\frac{C_n\tilde
\epsilon_{n+2}}{\lambda[n]^2
+C_n+\tilde \delta_n(\tilde
\delta_n-f_{n-1}(\lambda))+\tilde
\epsilon_n(\tilde \epsilon_n-e_{n-2}
(\lambda)\!-\!f_{n-1}(\lambda)(\tilde
\delta_{n-1}-f_{n-2} (\lambda)))
C_{n-1}^{-1}},
\end{equation}
with the initial conditions
\begin{equation} \label{4:10}
e_0(\lambda)=\tilde \epsilon_2,
\qquad e_1(\lambda)=\frac{C_1
\tilde \epsilon_{3}}{\lambda +
C_1}.
\end{equation}
\ec
\bdm
Using \eqref{4:4} and \eqref{4:7},
we get
$$\begin{array}{r@{}l}
{\bf d}^2_n= & \left\langle {Q^{(\lambda) }_n(z)},
{Q^{(\lambda)}_n(z)}\right\rangle_S=\left\langle
{Q^{(\lambda)}_n(z)},{p^{\{-1\}}_n(z)}\right\rangle_S
=d_n^2+\lambda\left\langle {{\bf u}}, {({\mathscr D}
Q^{(\lambda)}_n)(z)\,({\mathscr D}p^{\{-1\}}_n)(z)}
\right\rangle\\[4mm]
= & \displaystyle d^2_n+\lambda \left\langle {{\bf u}},
{({\mathscr D} Q^{(\lambda)}_n)(z)\,
\left([n]p_{n-1}(z)-\tilde\delta_{n}
({\mathscr D}p_{n-1})(z)-\tilde\epsilon_{n}
({\mathscr D}p_{n-2})(z)\right)}\right\rangle\\[4mm] =
& \displaystyle d^2_n+ \lambda [n]^2
d^2_{n\!-\!1}+\tilde \delta_n\left\langle {{\bf u}},
{Q^{(\lambda)}_n(z)\, p_{n-1}(z)}\right\rangle+\tilde
\epsilon_n \left\langle {{\bf u}},  {Q^{(\lambda)}_n
(z)\, p_{n-2}(z)}\right\rangle\\[3mm] = & d^2_n+ \lambda
[n]^2 d^2_{n\!-\!1}+\tilde \delta_n(\tilde \delta_n
\!-\!f_{n\!-\!1}(\lambda))d^2_{n\!-\!1} +\tilde \epsilon_n (\tilde
\epsilon_n\!-\! e_{n\!-\!2}(\lambda)\!-\!f_{n\!-\!1}(\lambda)
(\tilde \delta_{n\!-\!1}\!-\!f_{n\!-\!2}(\lambda))) d^2_{n\!-\!2},
\end{array} $$
 which proves \eqref{4:8}. Finally, \eqref{4:9} follows from
\eqref{4:8} using \eqref{3:6}, and the initial conditions
\eqref{4:10} can be deduced from \eqref{4:5}. \edm


The inner product  \eqref{4:1} does not satisfy the Hankel
property $\langle {zp},{r}\rangle_S=\langle {p},{zr} \rangle_S$,
and therefore the polynomial sequence $(Q_n^{(\lambda)})$ does not
fulfill in general a three--term recurrence relation. However, we
find an operator $\mathscr J$ which is symmetric with respect to
the inner product \eqref{4:1}, that is, $\langle {{\mathscr
J}p},{r}\rangle_S=\langle {p}, {{\mathscr J}r}\rangle_S,\, $ where
$ p,\, r\in \XX P.$ Thus, we generalize some results in this
direction obtained for the continuous  and discrete cases (see the
surveys \cite{maalre}, \cite{mamo}, \cite{fink1}, \cite{fink2},
\cite{mei1}, and the references therein, and also \cite{mayte}).
\bp \label{prop5.1} Let $\mathscr J$ be the linear operator
\begin{equation} \label{5:1} \mathscr J:=({\mathscr E}^-
\widetilde \phi)\mathscr I+\frac{\lambda}{q} ({\mathscr
D}^*\widetilde \phi -\widetilde \psi) {\mathscr D}^* -\lambda
({\mathscr E}^-\widetilde \phi){ \mathscr D}\,{\mathscr D}^*,
\end{equation}
where $\mathscr I$ is the identity operator. Then,
\begin{equation} \label{5:2}
\langle {({\mathscr E}^-\widetilde \phi) p},{r}\rangle_S=
\langle {{\bf u}}, {p\,\mathscr Jr}\rangle,\quad p,\, r \in \XX P.
\end{equation}
\ep \bdm According to Theorem \ref{the3.0} the linear functional
${\bf u}$ satisfies the distributional equation ${\mathscr D}^*
(\widetilde \phi{\bf u})=\widetilde \psi {\bf u}$. From
(\ref{4:1}), we have
$$
\langle {({\mathscr E}^-\widetilde \phi) p}, {r}\rangle_S=
\langle {{\bf u}}, {({\mathscr E}^-\widetilde \phi) p r}\rangle+
\lambda \langle {{\bf u}},{{\mathscr D}\big(({\mathscr E}^-
\widetilde \phi) p\big) \, {\mathscr D} r}\rangle.
$$
Now,  using relation \eqref{2:9} with $\pi={\mathscr D}^*(r)$ we
get
$$
{\mathscr D}\big(({\mathscr E }^-\widetilde \phi) p \, {\mathscr
D}^* r\big) = {\mathscr D} \big(({\mathscr E }^-\widetilde \phi)
p\big) \, {\mathscr D} r+ \big(({\mathscr E }^-\widetilde \phi)
p\big) {\mathscr D}\big({\mathscr D}^* r),
$$
and so,
$$
\langle {({\mathscr E}^-\widetilde \phi)p}, {r}\rangle_S= \langle
{{\bf u}}, {p\,\big(({\mathscr E}^-\widetilde \phi)r-\lambda
({\mathscr E}^-\widetilde \phi){\mathscr D}({\mathscr D}^*r)
\big)}\rangle+\lambda \langle {{\bf u}},{{\mathscr D}
\big(({\mathscr E}^-\widetilde \phi)p {\mathscr D}^*r)}\rangle.
$$
On the other hand,  if we apply the relation \eqref{2:8} to the
distributional equation with $\pi=\widetilde \phi$
we get
$$
\widetilde \psi {\bf u}=({\mathscr E}^-\widetilde \phi){\mathscr
D}^*{\bf u}+({\mathscr D}^*\widetilde \phi) {\bf u},
$$
and therefore
$$
\left \langle {({\mathscr E}^-\widetilde \phi)p}, {r}
\right\rangle_S= \displaystyle \left\langle {{\bf u}},
{p\, \big(({\mathscr E}^- \widetilde \phi) r-\lambda
({\mathscr E}^-\widetilde \phi) {\mathscr D}({\mathscr
D}^*r)\big)}\right\rangle -\frac{\lambda}{q} \langle{{\bf
u}} ,{p\, (\widetilde \psi-{\mathscr D}^*\widetilde \phi)
{\mathscr D}^*r}\rangle=\langle {{\bf u}} ,{p{\mathscr J}
r}\rangle.
$$
\edm \bc \label{cor5.1} The following identity holds:
$$
\langle {(\widetilde \psi-{\mathscr D}^*\widetilde \phi) p},
{r}\rangle_S= \langle {{\mathscr D}^*{\bf u}}, {p\, \mathscr
Jr}\rangle, \quad p,\, r\in \XX P.
$$
\ec \bdm Using the same properties as in the proof  of the
previous result, we get
$$\begin{array}{rl}
\langle {(\widetilde \psi-{\mathscr D}^*\widetilde \phi) p},
{r}\rangle_S=& \langle {{\bf u}}, {(\widetilde \psi-
{\mathscr D}^*\widetilde \phi) pr}\rangle+\lambda
\left\langle {{\bf u}},{{\mathscr D}\big(
(\widetilde \psi-{\mathscr D}^*\widetilde \phi) p\big)
\, {\mathscr D}r}\right\rangle\\[3mm] = &
\displaystyle \left\langle{{\mathscr D}^*{\bf u}}, {p
\left(({\mathscr E}^-
\widetilde \phi)r-\frac{\lambda}{q}(\widetilde \psi- {\mathscr
D}^*\widetilde \phi) {\mathscr D}^*r\right)}\right\rangle -\lambda
\langle {{\mathscr D}^*{\bf u}},{({\mathscr E}^-\widetilde \phi) p
{\mathscr D}({\mathscr D}^*r)}\rangle \\[3mm] = &
\displaystyle \langle {{\mathscr D}^*{\bf u}}, {p{\mathscr J} r}
\rangle.
\end{array} $$ \edm

In the next result we prove that the linear functional ${\mathscr
J}$ plays for the inner product \eqref{4:1} a role equivalent to
the one played by the multiplication operator by $z$ for standard
inner products (see (\ref{omi})).

\bt \label{the5.1} The linear operator ${\mathscr J}$ defined in
\eqref{5:1} satisfies
\begin{equation} \nonumber \label{5:3}
\langle {{\mathscr J}p},{r}\rangle_S=\langle {p},
{{\mathscr J}r}\rangle_S\qquad p,\, r\in \XX P.
\end{equation}
\et \bdm Applying Proposition \ref{prop5.1} and Corollary
\ref{cor5.1} together with (\ref{2:9}) and  the fact that
${\mathscr E}^+\, {\mathscr D}^*={\mathscr D}$, then  we get
$$
\begin{array}{rl}
\langle {{\mathscr J}p},{r}\rangle_S= &\displaystyle
\langle {({\mathscr E}^- \widetilde \phi) p},{r}\rangle_S
-\frac{\lambda}{q} \langle (\widetilde \psi-{\mathscr D}^*
\widetilde \phi)({\mathscr D}^*p),r \rangle_S-
\lambda \langle({\mathscr E}^-\widetilde \phi)
{\mathscr D}({\mathscr D}^*p),r\rangle_S
\\[3mm] =& \displaystyle \langle {{\bf u}},{p
{\mathscr J}r}\rangle-\frac{\lambda}{q} \langle {{\mathscr
D}^*{\bf u}},{({\mathscr D}^*p) {\mathscr J}r}\rangle -\lambda
\langle {{\bf u}},{({\mathscr D}({\mathscr
D}^*p)) {\mathscr J}r}\rangle\\[3mm]= &
\displaystyle \langle {{\bf u}},{p {\mathscr J}
r}\rangle+\lambda\langle {{\bf u}},{({\mathscr
D}p) ({\mathscr D} {\mathscr J}r)}\rangle
\\[3mm] = & \langle {p}, {{\mathscr J}r}\rangle_S.
\end{array}
$$
\edm

\bn \label{rem5.1} Taking into account the proof of the last
results, if $\gcd({\mathscr E}^-\widetilde \phi, {\mathscr D}^*
\widetilde \phi- \widetilde \psi)=d>1$, then we can consider
$({\mathscr E}^-\widetilde \phi)/d$ and $({\mathscr D}^*
\widetilde \phi-\widetilde \psi)/d$ in the definition (\ref{5:1})
of the operator $\mathscr J$ and these results also hold. \en

To end this section, we give some algebraic
differential/difference results. One of them allows us to give an
expression for the polynomial $({\mathscr E}^-\widetilde \phi)(z)
p_n(z)$ in terms of a finite number of ${\mathscr D}$--Sobolev
orthogonal polynomials. Another one gives the \textit{second--order
${\mathscr D}$--equation} (differential or difference or
$q$--difference equation) that  is satisfied by the polynomials
$(Q^{(\lambda)}_n)$.

\bc\label{prop5.2} The following relations hold
\begin{eqnarray*}  ({\mathscr
E}^-\widetilde \phi)(z) p_n(z)&= \sum_{\nu=n-H}^{\nu=n+\deg
\widetilde \phi} \mu_{n,\nu} Q^{(\lambda)}_\nu(z), \quad n\ge  H,
\\
{\mathscr J}Q^{(\lambda)}_n(z)&=\sum_{\nu=n-\deg
\widetilde\phi}^{n+H} \vartheta_{n,\nu}p_\nu(z),\quad n\ge \deg
\widetilde \phi,
\\
{\mathscr J}Q^{(\lambda)}_n(z)&= \sum_{\nu=n-H}^{n+H}\varpi_{n,
\nu}Q^{(\lambda)}_\nu(z), \quad n\ge H,
\end{eqnarray*}
where $H:=\max\{\deg \widetilde  \psi-1, \deg \widetilde \phi\}.$
 \ec \bdm To prove the first relation, we write
$$
({\mathscr E}^-\widetilde \phi)(z) p_n(z)=
\sum_{\nu=0}^{n+\deg \widetilde \phi}
\mu_{n,\nu} Q^{(\lambda)}_\nu(z),
$$
where $\mu_{n,\nu}{\bf d}^2_\nu=\langle {({\mathscr E}^-\widetilde
\phi) p_n} ,{Q^{( \lambda)}_\nu}\rangle_S=\langle {{\bf u}}, {p_n
{\mathscr J} Q^{(\lambda)}_\nu}\rangle$ which, using Proposition
\ref{prop5.1}, vanishes when $\nu+H<n$ since $\deg {\mathscr
J}\pi=H+\deg \pi$.

To prove the second relation it is enough to take (\ref{5:2}) into
account to get
$$
{\mathscr J}Q^{(\lambda)}_n(z)=\sum_{\nu=0}^{n+H}
\vartheta_{n,\nu}p_\nu(z),
$$
where $\vartheta_{n,\nu}d^2_\nu=\langle {{\bf u}}, {p_\nu
{\mathscr J}Q^{(\lambda)}_n}\rangle =\langle {({\mathscr
E}^-{\widetilde \phi}) p_\nu},{Q^{(\lambda)}_n}\rangle_S$ which
vanishes when $\nu+\deg \widetilde\phi<n$.

Finally,  using Theorem \ref{the5.1} we obtain
$$
{\mathscr J}Q^{(\lambda)}_n(z)= \sum_{\nu=0}^{n+H}\varpi_{n,\nu}
Q^{(\lambda)}_\nu(z),
$$
where $\varpi_{n,\nu}{\bf d}^2_\nu= \langle {{\mathscr
J}Q^{(\lambda)}_n}, {Q^{(\lambda)}_\nu}\rangle_S=\langle
Q^{(\lambda)}_n, {\mathscr J}Q^{(\lambda)}_\nu\rangle_S$ which
vanishes when $\nu+H<n, $ and thus the third relation is proved. \edm

\section{The examples}
We illustrate the results obtained in this paper with several
examples covering continuous, discrete and $q$--Hahn cases. Other
examples appearing in the literature are also included in this
general approach and their properties can be deduced from the
results introduced here. We also remark that the nonstandard inner
product associated with the $q$--Freud linear functional and
another one related to a $1$--singular semiclassical functional
are new and they have not been considered before.

\subsection{A continuous example:
the Jacobi--Sobolev Polynomials} The family of monic Jacobi
orthogonal polynomials $(P_n^{\alpha,\beta})$ is at the top of the
continuous classical polynomials in the Askey scheme and they can
be written as (see, for example, \cite{kost})
$$
P_n^{\alpha,\beta}(x)=\frac{2^n(\alpha+1)_n}{(\alpha+
\beta+n+1)_n}\,{}_2F_1\left.\left(\begin{array}{c}-n, \, \alpha
+\beta+n+1\\ \alpha+1\end{array} \right|\frac {1-x}2\right),
\qquad \alpha,\, \beta >-1.
$$
In fact, this family satisfies the following orthogonality
property:
$$
\left\langle {{\bf u}^{\alpha,\beta}},
{P^{\alpha,\beta}_n P^{\alpha,
\beta}_m}\right\rangle=\frac{2^{2n+\alpha+
\beta+1}\Gamma(n+1)\Gamma(\alpha+
n+1)\Gamma(\beta+n+1)\Gamma(\alpha+
\beta+n+1)}{(\alpha+\beta+2n+1)(\Gamma(
\alpha+\beta+2n+1))^2}\delta_{n,m},
$$
where the linear functional
${\bf u}^{\alpha,\beta}$ has
the following integral representation:
$$
\left\langle {{\bf u}^{\alpha,\beta}},P\right\rangle=
\int_{-1}^{1} P(x)(1-x)^\alpha
(1+x)^\beta dx,
$$
and satisfies the distributional equation:
\begin{equation} \nonumber
\mathscr{D}((1-x^2){\bf u}^{\alpha,
\beta})=(\beta-\alpha-x(\alpha+ \beta+2)){\bf u}^{\alpha,\beta},
\end{equation}
where, as we pointed out in Section 2, $\mathscr {D}=\mathscr
{D}^*=\frac d{dx}, \,$ so $t=2$, $p=1$, $\sigma=0$, and $H=2$.

Then, we can consider the Sobolev inner product defined by
\begin{equation}
\nonumber \langle f,g\rangle_S=\int_{-1}^1
f(x)g(x)(1-x)^\alpha(1+x)^\beta dx+ \lambda \int_{-1}^1
f'(x)g'(x)(1-x)^\alpha(1+x)^\beta dx,
\end{equation}
where $\alpha, \, \beta>-1$, and $\lambda\ge0$. This nonstandard
inner product has been considered in a lot of articles, and we
refer the reader to the surveys  mentioned several times during the
paper for more details. We  denote by $(Q_n^{\alpha,\beta})$ the
sequence of monic polynomials orthogonal with respect to the inner
product $(f,g)_S$, which are called  monic {\sf Jacobi--Sobolev
orthogonal polynomials}.

The results obtained in this paper allow us to recover some
relations between Jacobi and Jacobi--Sobolev orthogonal
polynomials:
 \begin{equation} \nonumber
P_n^{\alpha-1,\beta-1}(x)=Q_n^{\alpha,\beta}(x)+
\theta^{\alpha,\beta}_{n,n-1}Q_{n-1}^{\alpha,\beta}(x)+
\theta^{\alpha,\beta}_{n,n-2}Q_{n-2}^{\alpha,\beta}(x),
\end{equation}
\begin{equation} \nonumber
(x^2-1)P_n^{\alpha,\beta}(x)=Q^{\alpha,
\beta}_{n+2}(x)+\sum_{\nu=n-2}^{n+1}\mu^{\alpha,\beta}_{n,\nu}
Q^{\alpha, \beta}_{\nu}(x).
\end{equation}
Moreover, according to (\ref{5:1}) we can define  the linear
functional  \begin{equation} \nonumber {\mathscr
J}^{\alpha,\beta}=(1-x^2){\mathscr I}+\lambda(\alpha-
\beta+x(\alpha+\beta))\frac d{dx}-\lambda (1-x^2) \frac
{d^2}{dx^2}, \end{equation} and Corollary  \ref{prop5.2} yields
\begin{equation} \nonumber -{\mathscr J}^{\alpha,\beta}
Q_n^{\alpha,\beta}(x)=P_{n+2}^
{\alpha,\beta}(x)+\sum_{\nu=n-2}^{n+1}\vartheta^{\alpha,\beta}_{n,\nu}
P_{\nu}^{\alpha, \beta}(x),
\end{equation}
\begin{equation} \nonumber -{\mathscr
J}^{\alpha,\beta}Q_n^{\alpha,\beta}(x)=Q_{n+
2}^{\alpha,\beta}(x)+\sum_{\nu=n-2}^{n+1}\varpi^{\alpha,\beta}_{n,\nu}Q_{\nu}^{
\alpha,\beta}(x).
\end{equation}
Observe that the minus signs appear due to the factor
$(1-x^2){\mathscr I}$ in  $\mathscr J$.
\subsection{A discrete example: the $\Delta$--Meixner--Sobolev
polynomials}
Monic Meixner orthogonal  polynomials can be written as (see, for
example, \cite{kost})
$$
M_n(x;\beta,c)=\frac{(\beta)_nc^n}{(c-1)^n}\,
{}_2F_1\left.\left(\begin{array}{c}-n, \, -x\\ \beta\end{array}
\right| 1-\frac{1}{c} \right), \qquad \beta >0,\ 0<c<1.
$$
In fact, this family satisfies the following orthogonality
property:
$$
\left\langle {{\bf u}^{M}},{M_n M_m}\right\rangle=
\frac{(\beta)_nc^n n!}{(1-c)^{
2n+\beta}}\delta_{n,m},
$$
where the linear functional
${\bf u}^{M}$ has the following
representation:
$$
\left\langle {{\bf u}^{M}},P\right\rangle=\sum_{x=0}^\infty
P(x) \,\frac{(\beta)_x}{\Gamma(x+1)}\,c^x ,
$$
and it satisfies the distributional equations
\begin{equation} \nonumber
\label{7:13} \mathscr{D}^*((x+\beta){\bf u}^{M})=
\left(x\Big(1-\frac 1c\Big)+\beta \right) {\bf u}^{M},
\end{equation} \begin{equation}
\nonumber \label{7:13.1} \mathscr{D}(x{\bf u}^{M})=
\left(x(c-1)+\beta c\right) {\bf u}^{M},
\end{equation}
where $\mathscr {D}=\Delta$ and $\mathscr{ D}^*=\nabla$,
so $t=1$, $p=1$, $\sigma=0$, and $H=1$.

Now let us consider the nonstandard inner product defined by
\begin{equation}
\label{7:17} \langle f,g\rangle_S=\sum_{x=0}^\infty
f(x)g(x)\,\frac{(\beta)_x} {\Gamma(x+1)}\,c^x +\lambda \sum_{x=
0}^\infty (\Delta f)(x)(\Delta g)(x)\,
\frac{(\beta)_x}{\Gamma(x+1)}\,c^x,
\end{equation}
where $\beta>0$, $0<c<1$, and $\lambda>0$. This inner product is
known in the literature as the $\Delta$--Sobolev inner product.  We
denote by $(Q_n^{(\lambda)}(x;\beta,c))$ the  sequence of monic
polynomials orthogonal with respect to \eqref{7:17}, which are
called  monic {\sf $\Delta$--Meixner--Sobolev orthogonal
polynomials}.

Like in the Jacobi case, the results obtained in this paper allow
us to recover some relations between the families of orthogonal
polynomials $(M_n(x;\beta,c))$ and $(Q_n^{(\lambda)}(x;\beta,c))$:
\begin{equation} \nonumber \label{7:18}
M_n(x;\beta-1,c)=Q_n^{(\lambda)}(x;\beta,c)+
f^M_n(\lambda;\beta,c)Q_{n-1}^{(\lambda)} (x;\beta,c),
\end{equation}
\begin{equation}
\nonumber \label{7:21} (x-1)M_n(x;\beta,c)=Q^{(\lambda)}_{n+1}
(x;\beta,c)+\mu^M_{n,n}Q^{(\lambda)}_{n}
(x;\beta,c)+\mu^M_{n,n-1}Q^{\alpha, \beta}_{n-1}(x;\beta,c).
\end{equation}

Moreover we define the operator ${\mathscr J}^M$ as
\begin{equation}  \nonumber \label{7:19}
{\mathscr J}^M=(x-1){\mathscr I}+
\lambda\left(1+\beta c-x(c-1)\right)\nabla-\lambda (x-1)\Delta
\nabla,
\end{equation}
and applying the results of Section 3 we obtain
\begin{equation} \nonumber \label{7:22} {\mathscr J}^M
Q_n^{(\lambda)}(x;\beta,c)=M_{n+1}
(x;\beta,c)+\vartheta^M_{n,n}M_n(x;\beta,c)
+\vartheta^M_{n,n-1}M_{n-1}(x;\beta,c), \end{equation}
\begin{equation} \nonumber
\label{7:23} {\mathscr J}^M Q_n^{(\lambda)}(x;\beta,c)=Q_{n+
1}^{(\lambda)}(x;\beta,c)+\varpi^M_{n,n}
Q_n^{(\lambda)}(x;\beta,c)+\varpi^M_{n,n-1}
Q_{n-2}^{(\lambda)}(x;\beta,c). \end{equation}

\subsection{A $q$ example: $q$--Freud type polynomials}
The family of monic $q$--Freud polynomials, $(P_n)$, satisfies the
relation \cite{coma1}
\begin{equation} \nonumber \label{7:25} ({\mathscr
D}P_n)(x(s))=[n]P_{n-1} (x(s))+a_n P_{n-3}(x(s)), \qquad n\ge
0,\end{equation}
where $x(s)=q^s$, with $0<q<1$, ${\mathscr D}={\mathscr
D}_q$, $P_{-1}\equiv 0$, $P_{0}\equiv 1$, and $P_1(x)=x$.

So by Theorem \ref{lem3.1} we get $\phi(x)=1$, $\sigma=2$, and
$t=0$; hence $(P_n)$ is orthogonal with respect to the linear
functional ${\bf u}^{qF}$ of class 2  which fulfills the
distributional equations
\begin{equation} \label{7:27} {\mathscr D}({\bf
u}^{qF})=\psi {\bf u}^{qF}, \qquad \deg \psi=3, \end{equation}
\begin{equation} \nonumber
\label{7:27.1} {\mathscr D}^*((1+x(q-1)\psi){\bf u}^{qF})=q\psi
{\bf u}^{qF}.
\end{equation}
Furthermore, these polynomials are symmetric (see \cite{coma1})
and satisfy the three--term recurrence relation
$$
xP_n=P_{n+1}+c_nP_{n-1},\qquad n\ge 1. $$

In fact, a straightforward computation shows $a_n=K(q)
q^{-n}c_nc_{n-1}c_{n-2}$, and the sequence $(c_n)$ satisfies the
nonlinear recurrence relation
$$
q[n]c_{n-1}+K(q)q^{-n+1}c_nc_{n-1}c_{n-2}=[n-1]c_n+
K(q)q^{-n-1}c_{n+1}c_nc_{n-1},\quad n\ge 2,
$$
with initial condition $c_0=0$, and for the sake of simplicity
we choose the initial conditions $c_1$ and $c_2$ in such a way
that
$$
c_1^2+c_1c_2=1,
$$
and we have $\lim_{q\to 1^{-}}K(q)=4$. Therefore,
$\psi(x)=-K(q)q^{-3}x^3$ and the linear functional ${\bf u}^{qF}$
has the following integral representation:
\begin{equation} \nonumber
\label{7:28} \left \langle {{\bf u}^{qF}}, P\right\rangle=
\int_{-1}^{1} P(x)\,\frac{1}{((q-1)K(q)q^{-3}q^{4x};q^4)_{\infty}}
d_q(x),
\end{equation}
where
$$
\int_{-1}^{1} f(x) d_q(x)=(1-q)\sum_{k=0}^\infty f(q^k)q^k+
(1-q)\sum_{k=0}^\infty f(-q^k)q^k.
$$
In fact, this family satisfies the following orthogonality
property:
$$
\left\langle {{\bf u}^{qF}},{P_n P_m}\right\rangle=
2c_1c_2\cdots c_n\delta_{n,m}.
$$

Now, to illustrate our results we can introduce the nonstandard
inner product defined by
 \begin{equation} \label{7:32}
\begin{split} \langle f,g\rangle_S&=\int_{-1}^{1}
 f(x)g(x)\frac{1}{((q-1)K(q)q^{-3}q^{4x};q^4)_{\infty}} d_q(x)\\
 &+ \lambda \int_{-1}^{1} ({\mathscr D}_q f)(x)({\mathscr D}_q
f)(x)\frac{1}{((q-1)K(q)q^{-3}q^{4x};q^4)_{\infty}}
d_q(x).\end{split}
\end{equation}
We  denote by $(Q_n^{qF})$ the sequence of monic polynomials
orthogonal with respect to the inner product \eqref{7:32}, which
we  call the monic {\sf $q$--Freud--Sobolev orthogonal
polynomials}.

The theory developed in the previous sections  allows us to link
the $q$--Freud--Sobolev polynomials with the $q$--Freud
polynomials. Taking into account the distributional equation
\eqref{7:27} and Proposition \ref{prop5.1} we define
\begin{equation} \nonumber {\mathscr
J}^{qF}=(1-(q-1)K(q)q^{-7}x^4){\mathscr I}+\frac{\lambda}{q}
K(q)q^{-6}x^3{\mathscr D}_{1/q}-\lambda
(1-(q-1)K(q)q^{-7}x^4){\mathscr D}_{q}{\mathscr D}_{1/q}.
\end{equation}

In this case $H= \deg \widetilde \phi=4$ and the results in
Section 3 can be rewritten as \begin{equation} \nonumber
(1-(q-1)K(q)q^{-7}x^4)P_n(x)=\sum_{\nu=n-4}^{n+4}
\mu^{qF}_{n,\nu}Q^{qF}_{\nu} (x), \end{equation}
\begin{equation} \nonumber {\mathscr J}^{qF} Q_n^{qF}(x)=
\sum_{\nu=n-4}^{n+4}\vartheta^{qF}_{n,\nu}P_{\nu}(x),
\end{equation}
\begin{equation}
\nonumber {\mathscr J}^{qF} Q_n^{qF}(x)=
\sum_{\nu=n-4}^{n+4}\varpi^{qF}_{n,\nu}Q_{\nu}^{qF}(x).
\end{equation}

\subsection{A $1$--singular semiclassical
polynomials of class 1} The family considered for this example was
studied by J.C. Medem in \cite{med2}. Such a family of monic
polynomials, $(S_n)$, which we  call Medem polynomials, is
orthogonal with respect to the linear functional ${\bf w}$ which
satisfies the distributional equation
\begin{equation}  \label{7:40}
{\mathscr D}(x^3 {\bf w})=(-x^2+4){\bf w}, \end{equation} where
${\mathscr D}={\mathscr D}^*=\dfrac d{dx}$, $t=3$, $p=2$, and so
$\sigma=1$, with initial condition $({\bf w})_1=\langle {{\bf
w}},x\rangle=0$.

Moreover a straightforward computation shows that $a_{\sigma+1}+
b_{\sigma+2}=0$, i.e., ${\bf w}$ is $1$--singular. Indeed, ${\bf
w}$ is the symmetrized linear functional associated with the
linear functional ${\bf b}^{(-\frac 52)}$ (see \cite{ch} Chapter
1, Sections 8 and 9), i.e.,  $({\bf w})_{2n+1}=0$ and $({\bf
w})_{2n}= ({\bf b}^{(-\frac 52)})_n$, for any $n\ge 0$.

Notice that the linear functional ${\bf b}^{(\alpha)}$ has the
following integral representation:
$$
\left\langle {{\bf b}^{(\alpha)}}, P\right\rangle=\frac 1{2\pi i}
\int_{\XX T}P(z)z^\alpha e^{-\frac 2z}dz,\qquad \alpha>-2,
$$
and thus, after a straightforward calculation, we get that the linear
functional ${\bf w}$ has the following integral representation:
\begin{equation} \nonumber \label{7:41} \langle {{\bf
w}},P\rangle=\frac 1{2\pi i}\int_{\XX T} P(z) z^{-4}
e^{-\frac 2{z^{2}}}dz,
\end{equation}
where $\XX T:=\{z\in \XX C: |z|=1\}$ is the unit circle.

This family can be written in terms of the Bessel
polynomials as follows \cite{med2}:
\begin{align*}
S_{2n}(x)=& \displaystyle \frac{2^n}{(n-\frac 32)_n}
B^{(-\frac 52)}_n(x^2)=r_{2n}x^{5-n}e^{\frac 2{x^2}}
({\mathscr D}^*)^n\left(x^{4n-5}e^{-\frac 2{x^2}}\right),
\\ S_{2n+1}(x)=& \displaystyle
\frac{2^n}{(n-\frac 12)_n}xB^{(-\frac 32)}_n(x^2)
=r_{2n+1}x^{3-n}e^{\frac 2{x^2}}
({\mathscr D}^*)^n\left(x^{4n-3}e^{-\frac 2{x^2}}\right),
\end{align*}

where $B_n^{(\alpha)}$ is the Bessel polynomial of degree $n$ with
parameter $\alpha$, and $r_n\ne 0$ are the corresponding
normalization coefficients for $n\ge 0$.

We can introduce the nonstandard inner product defined by
\begin{equation} \label{7:47}
(f,g)_S=\frac 1{2\pi i}\int_{\XX T} f(z)\overline
{g(z)}z^{-4}e^{-\frac 2{z^{2}}}dz+\frac{\lambda}{2\pi i}\int_{\XX T}
f'(z)\overline{g'(z)}z^{-4}e^{-\frac2{z^{2}}}dz.
\end{equation}
We  denote by $(Q_n^{S})$ the sequence of monic polynomials
orthogonal with respect to the inner product \eqref{7:47}, which
we  call the monic {\sf Medem--Sobolev orthogonal polynomials}.

Taking into account the distributional equation \eqref{7:40}  and
Proposition \ref{prop5.1} we define
\begin{equation} \nonumber \label{7:34}
{\mathscr J}^{S}=x^3{\mathscr I}+\lambda
(4x^2-4)\frac{d}{dx}-\lambda x^3\frac{d^2}{dx^2}.
\end{equation}

In this case $H= \deg \widetilde \phi=3$ and the results in
Section 3 can be rewritten as

\begin{equation} \nonumber \label{7:36}
x^3S_n(x)=Q^{S}_{n+3}(x)+\sum_{\nu=n-3}^{n+2}\mu^{S}_{n,\nu}
Q^{S}_{\nu}(x), \end{equation} \begin{equation} \nonumber
\label{7:37} {\mathscr J}^{S} Q_n^{S}(x)= S_{n+3}(x)+
\sum_{\nu=n-3}^{n+2}\vartheta^{S}_{n,\nu}
S_{\nu}(x), \end{equation} \begin{equation} \nonumber
\label{7:38} {\mathscr J}^{S} Q_n^{S}(x)=Q^{S}_{n+3}(x)+
\sum_{\nu=n-3}^{n+2}\varpi^{S}_{n,\nu}
Q^{S}_{\nu}(x). \end{equation}


{\bf Acknowledgements}

We thank the anonymous referees for
the revision of the manuscript. Their comments and suggestions have
improved the presentation of the work. The first author was
partially supported by MICINN of Spain under Grant MTM2006-13000-C03-02.
The second author was partially supported by MICINN of Spain under
Grant MTM2008-06689-C02-01 and Junta de Andaluc\'ia
(FQM229 and excellence project P06-FQM-1735).

\newpage

{\large \textbf{Appendix. Proof of Theorem \ref{lem3.2}}}

We are going to prove the result for $k\ge 0$ since the case $k<0$
is totally analogous to this one. When $k=0$ we know that $\phi$
and $\psi$ are polynomials  of degree, at most $\sigma+2$ and
$\sigma+1$, respectively, where $\sigma=\sigma_0$ is the order of
$\bf u$ with respect to the pair $(\phi,\psi)$. Then, using
(\ref{2:7}) we get for any monic polynomial $\pi$,
$$
{\mathscr D}(\pi \phi {\bf u}_0)= ({\mathscr D}\pi)\phi {\bf u}_0
+(\mathscr{E}^+\pi){\mathscr D}(\phi {\bf u}_0) =\big(({\mathscr
D}\pi)\phi +({\mathscr E}^+\pi)\psi\big){\bf u}_0,
$$
where $\tilde{\pi}:=({\mathscr D}\pi)\phi +({\mathscr
E}^+\pi)\psi$ is a polynomial of degree, at most,
$\deg \pi+\sigma+1$. \\
Moreover, if the pair $(\phi,\psi)$ is admissible, then the
leading coefficient is
$$
q^{m}(a_{\sigma+1}+q^{-1}[m]^* b_{\sigma+2})\ne 0, \quad
\mathrm{with }\quad m=\deg \pi.
$$

If we assume the result holds  for every $0\le k<K$, let us prove
it for $k=K.$ Thus, the linear functionals ${\bf u}_k$ satisfy the
distributional equation ${\mathscr D} (\phi_k {\bf u}_k)=\psi_k
{\bf u}_k$ for $0\le k<K$, and then we can obtain
$$\begin{array}{rl}
\displaystyle {\mathscr D}(\phi_{K-1}{\bf u}_K)& \displaystyle =
({\mathscr D}\phi_{K-1}){\bf u}_K+ ({\mathscr
E}^+\phi_{K-1}){\mathscr D}{\bf u}_K
\\[2mm] & = ({\mathscr D}\phi_{K-1})
{\bf u}_K+q({\mathscr E}^+\phi_{K-1})\mathscr{E}^+
\left(\left(\mathscr{D}\phi_{K-1}\right) {\bf u}_{K-1}\right)
\\[2mm] & = ({\mathscr D}\phi_{K-1}){\bf u}_K+
q({\mathscr E}^+\phi_{K-1}){\mathscr E}^+
\big(\psi_{K-1} {\bf u}_{K-1}\big)\\[2mm] & =
({\mathscr D}\phi_{K-1}){\bf u}_K+ q({\mathscr E}^+
\psi_{K-1}){\mathscr E}^+\big(\phi_{K-1}{\bf u}_{K-1}
\big)\\[2mm] & = \big({\mathscr D}\phi_{K-1}+q({\mathscr
E}^+\psi_{K-1})\big){\bf u}_K=\widetilde \psi_K\, {\bf u}_K.
\end{array}
$$ Therefore  taking  into account Definition   \ref{def3.4} and Remark \ref{rem3.3}, there
exist polynomials $\phi_K, \psi_K$ satisfying
\begin{equation} \label{3:15.1} {\mathscr
D}(\phi_K {\bf u}_K)= \psi_K {\bf u}_K,
\end{equation}
and a monic polynomial $\xi_K$ such that $\phi_{K-1}=\xi_K
\phi_{K}$, and now, taking (\ref{3:15.1}) into account
 the polynomials $\psi_K$ and $\psi_{K-1}$ fulfill
the relation
\begin{equation} \label{3:15}
{\mathscr D}(\phi_{K-1})+q({\mathscr E}^+\psi_{K-1})= {\mathscr
D}(\xi_K)\phi_{K}+{\mathscr E}^+(\xi_K)\psi_{K}, \end{equation}
and we can check that $\psi_K$ is a polynomial of degree, at most,
$\sigma_{K}+1$. Furthermore, using \eqref{3:15} and working in the
same way as in the case $k=0$, we get
\begin{equation} \nonumber \label{3:16}
{\mathscr D}(\pi \phi_K {\bf u}_K)= \big(({\mathscr D}\pi)\phi_K
+({\mathscr E}^+\pi)\psi_K\big){\bf u}_K. \end{equation}  Clearly,
$\tilde{\pi}:=({\mathscr D}\pi)\phi_K +({\mathscr E}^+\pi)\psi_K$
is a polynomial of degree, at most, $m+\sigma_K+1$.

In fact, if the pair $(\phi_K, \psi_K)$ is admissible the leading
coefficient of $\widetilde \pi$ is
\begin{equation} \label{3:16.1}
q^{m}(a_{\sigma_K+1}+q^{-1}[m]^* b_{\sigma_K+2})\ne 0,
\end{equation}
where the coefficients $a_{\sigma_K+1}$ and $b_{\sigma_K+2}$ can
be obtained recursively
\begin{equation} \label{3:16.3}
q^{\sigma_{K-1}+2}\left(a_{\sigma_{K-1}+1}+q^{-1}
[\sigma_{K-1}+2]^*b_{\sigma_{K-1}+2}\right)=
q^{\sigma_{K-1}-\sigma_{K}}\left(a_{\sigma_{K}+1}
+q^{-1}[\sigma_{K-1}-\sigma_{K}]^*b_{\sigma_{K}+2}\right).
\end{equation}
Notice that in the continuous and discrete cases $q=1$ and
$[m]^*=m.$ Then, using the admissibility condition of Definition
\ref{def3.2}, the above expression, and (\ref{3:16.1}) in a
recursive way, we get that the leading coefficient is equal to
\begin{equation} \nonumber \label{3:16.2}
a_{\sigma_K+1}+m
b_{\sigma_K+2}=a_{\sigma+1}+\left(m+\sum_{\nu=1}^{k}(\sigma_\nu+2)\right)
b_{\sigma+2}\ne 0, \quad k\ge 0.
\end{equation}

Therefore in the $q$--Hahn case there are infinitely many values
of $q$ for which the expression (\ref{3:16.1}) is different from
zero for any  nonnegative integers $m$ and $k$. Hence the result
follows for $k\ge 0$.

Observe that by \eqref{3:16.3} if there exists an integer $k$ such
that the  pair $(\phi_k,\psi_k)$ is admissible, then the  pair
$(\phi_\ell,\psi_\ell)$ is admissible for every $\ell\ge k$.

\end{document}